\documentclass[a4paper,12pt]{article}

\usepackage{amsmath}
\usepackage{amsthm}
\usepackage{amssymb}  
\usepackage{latexsym} 
\usepackage[all]{xy}
\usepackage{comment}
\usepackage{rotating}
\usepackage{slashbox}
\usepackage{multirow}
\usepackage{rotating}
\usepackage{stmaryrd}

\newcommand{\Mathematica}{{\sf Mathematica}}

\newcommand{\Q}{{\mathbb Q}}

\newcommand{\Z}{{\mathbb Z}}

\newfont{\wncyr}{wncyr10 at 12pt}
\newfont{\wncyrten}{wncyr10 at 10pt}

\newenvironment{Proof}{\par\noindent{\sc Proof:}}%
                      {\hspace*{\fill}\nobreak$\Box$\par\medskip}
                       {\hspace*{\fill}\nobreak$\Box$\par\medskip}
\newenvironment{myitemize}
{\begin{itemize}
\setlength{\itemsep}{1pt}
\setlength{\parskip}{0pt}
\setlength{\parsep}{0pt}}
{\end{itemize}}

\newtheorem{Proposition}{Proposition}[section]
\newtheorem{Theorem}[Proposition]{Theorem}
\newtheorem{Lemma}[Proposition]{Lemma}
\newtheorem{Corollary}[Proposition]{Corollary}

\theoremstyle{definition}

\newtheorem{Remark}[Proposition]{Remark}

\newcounter{nootje}
\begin{document}
\normalsize
\title{Formal groups and invariant differentials of elliptic curves}
\author{Mohammad Sadek}
\date{}
\maketitle
\let\thefootnote\relax\footnote{Mathematics Subject Classification: 14H52}
\begin{abstract}{\footnotesize In this paper, we find a power series expansion of the invariant differential $\omega_E$ of an elliptic curve $E$ defined over $\Q$, where $E$ is described by certain families of Weierstrass equations. In addition, we introduce several congruence relations satisfied by the trace of the Frobenius endomorphism of $E$.}
\end{abstract}
\textbf{Keywords:} Elliptic curves; Formal groups; Invariant differentials.
\section{Introduction}

Let $E$ be an elliptic curve defined over $\Q$ described by the Weierstrass equation \[y^2+a_1xy+a_3y=x^3+a_2x^2+a_4x+a_6,a_i\in\Z.\] Choosing a local parameter $z=-x/y$ for $E$ at its origin $O_E$, one can associate to $E$ a power series $\displaystyle w(z)=\sum_{n=0}^{\infty}s_iz^i\in\Z[a_1,a_2,a_3,a_4,a_6]\llbracket z\rrbracket$. Consequently, there are Laurent series expansions for the coordinates $x(z)$ and $y(z)$, hence one obtains power series expansions for several arithmetic objects attached to $E$ including the invariant differential $\omega_E$, the formal logarithm $\log_E(z)=\displaystyle\int\omega_E(z)$, and the formal group associated to $E$ over $\Z$ given by $\displaystyle F(X,Y)=\log_E^{-1}\left(\log_E(X)+\log_E(Y)\right)$.

Honda in \cite{HondaFormal} found an interesting link between the $L$-series $\displaystyle L(s)=\sum_{n=1}^{\infty}c_nn^{-s}$ of an elliptic curve $E$ and the formal group associated to $E$. If one sets $\displaystyle g(x)=\sum_{n=1}^{\infty}n^{-1}c_nx^n$, then $G(X,Y)=g^{-1}(g(X)+g(Y))$ is a formal group over $\Z$; moreover, $G(X,Y)$ is isomorphic to the formal group law $F(X,Y)$ associated to $E$ over $\Z$.
The isomorphism between these formal groups is made explicit to produce Atkin-Swinnerton-Dyer congruence relations. These congruence relations connect the coefficients of the $L$-series and the coefficients of the power series expansion of the invariant differential $\omega_E$. The modularity of elliptic curves $E$ defined over $\Q$ implies the existence of a set of congruence relations between the coefficients of the modular form attached to $E$ and the coefficients of the power series expansion of $\omega_E$.

The above discussion indicates that explicit formulas for the coefficients of the power series of $\omega_E$ will provide us with information about formal groups, $L$-series and modular forms associated to elliptic curves. There are a few explicit descriptions for the power series expansion of $\omega_E$ of a certain elliptic curve $E$. The expansion of $\omega_E$ can be found in \cite{CosterLegendre} when $E$ has a rational $2$-torsion point, i.e., when $E$ is described by a Weierstrass equation of the form $\displaystyle y^2=x^3+a_2x^2+a_4x$. It is shown that $\displaystyle \omega_E=\sum_{n=0}^{\infty}P_n\left(a_2/\sqrt{\Delta}\right)\left(\sqrt{\Delta}\right)^n\,z^{2n}\,dz,$ where $P_n$ is the $n$-th Legendre polynomial and $\Delta=a_2^2-a_4$. This enables the authors to find explicit congruence relations satisfied by Legendre polynomials using Atkin-Swinnerton-Dyer congruences and sharper congruences using other techniques when $E$ has complex multiplication. Most recently, an explicit $t$-expansion of $\omega_E$ was found for the elliptic curve $E:y^2=4(x^3+Ax+B)$ where $t=-2x/y$, see \cite{Yasuda}.

More combinatorial quantities appear as coefficients of invariant differentials. In \cite{stienstraLfunctions}, the integers $\displaystyle (-1)^m\sum_{k}{m\choose k}^3$ turn out to be the coefficients of the holomorphic differential form of a model of a $K3$-surface. In \cite{BeukersStienstra}, Ap\'{e}ry Numbers, $\displaystyle \sum_k{n\choose k}^2{n+k\choose k}$, which were used in Ap\'{e}ry's proof of the irrationality of $\zeta(2)$ and $\zeta(3)$, appear as the coefficients of the holomorphic differential form of a model of a $K3$-surface.

In this paper, we write explicit formulas for the $z$-power series expansions of the invariant differentials of elliptic curves described by the Weierstrass equations $y^2+a_1xy+a_3y=x^3+a_2x^2+a_4x$ and $y^2+a_3y=x^3+a_6$. We find a power series solution to the functional equation satisfied by $w(z)$ and use this to write a $z$-power series expansion for $\omega_E$.
Several combinatorial numbers appear as coefficients of these power series. These numbers include
\begin{align*}
&\sum_{m=\lfloor n/2\rfloor}^{n}\sum_{k=0}^{\lfloor m/2\rfloor}\sum_{r=0}^{n-m-k}{m\choose k}{m-k\choose k}{m-2k\choose r}{k\choose n-m-k-r},{\textrm{ and}}\\
&\sum_{k=\lfloor(n-1)/2\rfloor}^{n-1}{n+k-1\choose k}{k\choose n-k-1},
\end{align*}
thus we have families of congruence relations satisfied by the coefficients of the modular forms of these elliptic curves, and an explicit description of the formal logarithms of the formal groups associated to them.


It is worth mentioning that although one can find for any Weierstrass equation the power series $w(z)$, it is not clear for the author how to use it to obtain simple formulas for the invariant differential when $a_i\ne0$ for every $i$. However, the two families of Weierstrass equations that we treat are broad enough to include elliptic curves with non-trivial rational points for example.

\section{Formal groups of elliptic curves}
Let $E$ be an elliptic curve defined over the rational field $\Q$. Assume that $E$ is described by the following Weierstrass equation \[y^2+a_1xy+a_3y=x^3+a_2x^2+a_4x+a_6,\;a_i\in\Z.\]
We furthermore assume that the Weierstrass equation is globally minimal.
Let $A$ be the local ring of functions defined at the origin $O_E$, and $\widehat{A}$ the completion of $A$ at its maximal ideal. Then $\widehat{A}$ is isomorphic to the power series ring $\Q\llbracket z\rrbracket$ where $z$ is a parameter at the origin. This is used to express the Weierstrass coordinates as formal power series in $z$. More explicitly, we set \[z=-\frac{x}{y},\textrm{ and }w=-\frac{1}{y}\]
Now $z$ is a parameter at the origin $(z,w)=(0,0)$. The Weierstrass equation for $E$ becomes
\[w=z^3+a_1zw+a_2z^2w+a_3w^2+a_4zw^2+a_6w^3\] substituting the equation into itself recursively we write $w$ as a power series in $z$
\[w(z)=\sum_{i\ge0}s_i z^i\in\Z[a_1,\ldots,a_6]\llbracket z\rrbracket\] where $s_0=s_1=s_2=0,\;s_3=1$, and for $n\ge 4$
\begin{eqnarray}
\label{eq1}
s_n&=&a_1s_{n-1}+a_2s_{n-2}+a_3\sum_{k+l=n}s_ks_l+a_4\sum_{k+l=n-1}s_ks_l+a_6\sum_{k+l+m=n}s_ks_ls_m.
\end{eqnarray}
The reader can consult \cite{AntoniaBluher} and Chapter IV in \cite{sil1}. Formal series expressions for $x$ and $y$ can be deduced as $\displaystyle x(z)=\frac{z}{w(z)}$ and $\displaystyle y(z)=\frac{-1}{w(z)}$. Therefore, the pair $(x(z),y(z))$ is a formal solution to the Weierstrass equation of $E$.

The invariant differential $\displaystyle\omega_E$ of $E$ can be expressed as a formal power series in $z$
\[\omega_E=\frac{dy}{3x^2+2a_2x+a_4-a_1y}=\sum_{n=1}^{\infty}b(n)z^{n-1}\,dz,\;b(n)\in\Z,\;b(1)=1.\]

Let $p$ be a prime of good reduction for $E$. The trace of the Frobenius endomorphism modulo $p$ is $t_p=1+p-\#E(\mathbb{F}_p)$. The following congruences are Atkin-Swinnerton-Dyer congruences modulo a prime of good reduction, see \cite{Hazewinkellecturesformal}.

\begin{Corollary}
\label{cor1}
If $E$ has good reduction mod $p$, then
\begin{myitemize}
\item[a)] $b(p)\equiv t_p$ mod $p$;
\item[b)] $b(np)\equiv b(n)\,b(p)$ mod $p$ if $p\nmid n$;
\item[c)] $b(np)-t_p\,b(n)+p\,b(n/p)\equiv 0\textrm{ mod }p^s\textrm{ if }n\equiv 0\textrm{ mod }p^{s-1},s\ge 1$.
\end{myitemize}
\end{Corollary}

\section{The elliptic curve $y^2+a_1xy+a_3y=x^3+a_2x^2+a_4x$}
In this section we consider elliptic curves over $\Q$ given by the Weierstrass equation $y^2+a_1xy+a_3y=x^3+a_2x^2+a_4x$ where $a_i\in \Z$. Any elliptic curve with a non-trivial torsion point can be described by such a Weierstrass equation. Set $z=-x/y$ to be a parameter at the origin, and $w=-1/y$. According to equation (\ref{eq1}), $\displaystyle w(z)=\sum_{n=0}^{\infty}s_nz^n$, where $s_0=s_1=s_2=0,s_3=1$ and
\[s_n=a_1s_{n-1}+a_2s_{n-2}+a_3\sum_{k+l=n}s_ks_l+a_4\sum_{k+l=n-1}s_ks_l.\]
 The generating function $w(z)$ of the sequence $\displaystyle (s_n)_{n=0}^{\infty}$ satisfies the following functional equation:
 \[w(z)=z^3+a_1zw(z)+a_2z^2 w(z)+a_3w(z)^2+a_4 z w(z)^2.\]
 The above equation is quadratic in $\omega(z)$. As a consequence, one has
 \begin{eqnarray}
 \label{eq2}
 w(z)=\frac{1-a_1z-a_2z^2-\sqrt{(1-a_1z-a_2z^2)^2-4z^3(a_3+a_4z)}}{2(a_3+a_4z)}.
 \end{eqnarray}
We chose the negative sign because the positive sign would force $w(z)$ to have a pole at $z=0$.
Recall that $\displaystyle x(z)=z/w(z)$ and $y(z)=-1/w(z).$ The invariant differential $\omega_E$ of $E$ is given by
\begin{eqnarray}
\label{eq:invariant}
\omega_E&=&\frac{dy}{3x^2+2a_2x+a_4-a_1y}=\frac{\frac{d}{dz}\,w}{3z^2+2a_2zw+a_4w^2+a_1w}\,dz.
\end{eqnarray}
 We use \Mathematica, \cite{Mathematica}, to substitute (\ref{eq2}) in the formula for $\omega_E$ to obtain
 \begin{eqnarray*}
 \omega_E &=&\frac{1}{\sqrt{1-2(a_1+a_2z)z+\left[(a_1+a_2z)^2-4(a_3+a_4z)z\right]z^2}}\,dz.
 \end{eqnarray*}

If $\phi(x)$ is a power series, then we write $[x^n]\phi(x)$ for the coefficient of $x^n$ in $\phi(x)$. Before we proceed to the main result of this section we recall that the generalized central trinomial polynomials $T_n(x,y)$ is given as follows
\begin{eqnarray*}
T_n(x,y)=[t^n](t^2+xt+y)^n=\sum_{k=0}^{\lfloor n/2\rfloor}{n\choose k}{n-k\choose k}x^{n-2k}y^k,
\end{eqnarray*}
where the generating function is given by
$\displaystyle
\frac{1}{\sqrt{1-2xt+(x^2-4y)t^2}}=\sum_{n=0}^{\infty}T_n(x,y)t^n.$
\begin{Theorem}
\label{thm:family1}
Let $E:y^2+a_1xy+a_3y=x^3+a_2x^2+a_4x,\;a_i\in\Z,$ be an elliptic curve over $\Q$. Let $z=-x/y$ be a local parameter at $O_E$. The $z$-power series expansion of $\omega_E$ is given by $\displaystyle\sum_{n=0}^{\infty}b(n+1)\,z^n\,dz$, where $b(n+1)$ is
\[\sum_{m=\lfloor n/2\rfloor}^{n}\sum_{k=0}^{\lfloor m/2\rfloor}\sum_{r=0}^{n-m-k}{m\choose k}{m-k\choose k}{m-2k\choose r}{k\choose n-m-k-r}a_1^{m-2k-r}a_2^ra_3^{2k-n+m+r}a_4^{n-m-k-r}\]
\end{Theorem}
\begin{Proof}
We have already shown that $\displaystyle \omega_E=\frac{1}{\sqrt{1-2(a_1+a_2z)z+\left[(a_1+a_2z)^2-4(a_3+a_4z)z\right]z^2}}\,dz$. We compare the formula for $\omega_E$ with the generating function of the generalized central trinomial polynomials. One concludes that
\begin{eqnarray*}
\omega_E&=&\sum_{m=0}^{\infty}T_m\left(a_1+a_2z,(a_3+a_4z)z\right)\,z^m\,dz.
\end{eqnarray*}
Now one has
\begin{eqnarray*}
&T_m&\textrm{\hskip-10pt}(a_1+a_2z,(a_3+a_4z)z)= \sum_{k=0}^{\lfloor m/2\rfloor}{m\choose k}{m-k\choose k}(a_1+a_2z)^{m-2k}(a_3+a_4z)^kz^k\\
&=&\sum_{k=0}^{\lfloor m/2\rfloor}{m\choose k}{m-k\choose k}\left(\sum_{i=0}^{m-2k}{m-2k\choose i}a_1^{m-2k-i}a_2^i\,z^i\right)\left(\sum_{j=0}^k{k\choose j}a_3^{k-j}a_4^jz^{j}\right)z^k\\
&=&\sum_{k=0}^{\lfloor m/2\rfloor}{m\choose k}{m-k\choose k}\sum_{s=0}^{m-k}\left(\sum_{r=0}^s{m-2k\choose r}{k\choose s-r}a_1^{m-2k-r}a_2^ra_3^{k-s+r}a_4^{s-r}\right)z^{s+k}.
\end{eqnarray*}
For the third equality we applied the formula for multiplying polynomials. Now,
\begin{eqnarray*}
&b(n+1)&=[z^n]\,\omega_E=[z^n]\sum_{m=0}^{\infty}T_m\left(a_1+a_2z,(a_3+a_4z)z\right)\,z^m\\
&=&[z^n]\sum_{m=0}^{\infty}\sum_{k=0}^{\lfloor m/2\rfloor}\sum_{s=0}^{m-k}\sum_{r=0}^s{m\choose k}{m-k\choose k}{m-2k\choose r}{k\choose s-r}a_1^{m-2k-r}a_2^ra_3^{k-s+r}a_4^{s-r}z^{s+k+m}.
\end{eqnarray*}
Therefore, $b(n+1)$ consists of the sum of the coefficients for which $s+k+m=n$. Since $0\le s\le m-k$, one has $n\le 2m$. Moreover, $m$ cannot exceed $n$. It follow that $b(n+1)$ is
\begin{eqnarray*}
\sum_{m=\lfloor n/2\rfloor}^{n}\sum_{k=0}^{\lfloor m/2\rfloor}\sum_{r=0}^{n-m-k}{m\choose k}{m-k\choose k}{m-2k\choose r}{k\choose n-m-k-r}a_1^{m-2k-r}a_2^ra_3^{2k-n+m+r}a_4^{n-m-k-r}.
\end{eqnarray*}
\end{Proof}

The following corollary presents the Atkin-Swinnerton-Dyer congruences satisfied by $b(n+1)$ modulo the primes of good reduction of $E$.
\begin{Corollary}
\label{cor:family1}
Let $E/\Q$ be an elliptic curve defined by $y^2+a_1xy+a_3y=x^3+a_2x^2+a_4x$. Let $t_p=1+p-\#E(\mathbb{F}_p)$ where $p$ is a prime of good reduction for $E$. Assume further that $\omega_E=\sum_{n=1}^{\infty}b(n)\,z^{n-1}\,dz$ where $z=-x/y$ and $b(n)$ is given in Theorem \ref{thm:family1}. Then
\begin{myitemize}
\item[a)] $b(p)\equiv t_p$ mod $p$;
\item[b)] $b(np)\equiv b(n)\,b(p)$ mod $p$ if $p\nmid n$;
\item[c)] $b(np)-t_p\,b(n)+p\,b(n/p)\equiv 0\textrm{ mod }p^s\textrm{ if }n\equiv 0\textrm{ mod }p^{s-1},s\ge 1$.
\end{myitemize}
\end{Corollary}
\begin{Proof}
This is Corollary \ref{cor1}.
\end{Proof}
Using Theorem \ref{thm:family1}, one has the following explicit description of the $z$-power series expansion of $\omega_E$ for any elliptic curve $E/\Q$ described by $E:y^2+a_3y=x^3$. In fact, one obtains
\begin{eqnarray*}
\omega_E=\frac{1}{\sqrt{1-4a_3z^3}}\,dz=\sum_{n=0}^{\infty}{2n\choose n}a_3^nz^{(3n+1)-1}\,dz.
\end{eqnarray*}
Observing that the discriminant of $E$ is $ \Delta_E=-27a_3^4$, one uses Corollary \ref{cor:family1} to obtain the following congruences.
\begin{Corollary}
\label{corsuper}
 Let $E/\Q$ be an elliptic curve defined by $y^2+ay=x^3$ and $p$ a prime such that $p\nmid 3a$. One has
 $\displaystyle t_p= 0$ if $p\equiv 2$ mod $3$.
If $p\equiv 1$ mod $3$, then
\begin{myitemize}
\item[a)] $\displaystyle t_p\equiv {2\left(\frac{p-1}{3}\right)\choose\frac{p-1}{3}}a^{(p-1)/3}$ mod $p$;
\item[b)] $\displaystyle {2\left(\frac{np-1}{3}\right)\choose\frac{np-1}{3}}\equiv{2\left(\frac{n-1}{3}\right)\choose\frac{n-1}{3}}{2\left(\frac{p-1}{3}\right)\choose \frac{p-1}{3}}$ mod $p$ if $n\equiv 1$ mod $3$.
\item[c)] $\displaystyle {2\left(\frac{np-1}{3}\right)\choose \frac{np-1}{3}}a^{(np-1)/3}-t_p{2\left(\frac{n-1}{3}\right)\choose \frac{n-1}{3}}a^{(n-1)/3}+p{2\left(\frac{n/p-1}{3}\right)\choose \frac{n/p-1}{3}}a^{(n/p-1)/3}\equiv0\textrm{ mod }p^{s},$ if $n\equiv0$ mod $ p^{s-1},\;s\ge1$.
\end{myitemize}
\end{Corollary}
As pointed out by the referee, the congruence in Corollary \ref{corsuper} b) holds modulo $p^2$, see for example \cite{Coster2} where $p$-adic techniques are used to produce such congruence relations modulo higher powers of $p$. Investigating similar congruence relations modulo higher powers of $p$ will be the subject of future work.
 
\section{The elliptic curve $y^2+a_3y=x^3+a_6$}
We consider the family of elliptic curves defined by $y^2+a_3y=x^3+a_6$, $a_i\in\Z$. We set $z=-x/y$ and $\displaystyle w(z)=-1/y=\sum_{n=0}^{\infty}s_nz^n$, where $s_0=s_1=s_2=0,s_3=1$, and $\displaystyle s_n=a_3\sum_{k+l=n}s_ks_l+a_6\sum_{k+l+m=n}s_ks_ls_m$, $n\ge4$. Now $w(z)$ satisfies the following functional equation
$\displaystyle
w(z)=z^3+a_3w(z)^2+a_6w(z)^3.$
The latter equation can be written as
\[v=z^3=w(z)\left(1-a_3w(z)-a_6w(z)^2\right)=\frac{w(z)}{1/\left(1-a_3w(z)-a_6w(z)^2\right)}.\]
In order to find the power series expansion of $\omega_E$, we will need the following lemma.
\begin{Lemma}[Lagrange Inversion Theorem]
\label{lem:lagrange}
Suppose $u = u(x)$ is a power series in $x$ satisfying $x = u/\phi(u)$ where $\phi(u)$ is a power series in $u$ with a nonzero constant term. Then we have
\[[x^n]u(x)=\frac{1}{n}[u^{n-1}]\phi^{n}(u).\]
\end{Lemma}
Using Lemma \ref{lem:lagrange}, one obtains
\begin{eqnarray*}
[v^n]w(v)&=&\frac{1}{n}\,[w^{n-1}] \left(\frac{1}{1-a_3w-a_6w^2}\right)^n
=\frac{1}{n}\,[w^{n-1}]\sum_{k=0}^{\infty}(-1)^k{-n\choose k}(a_3+a_6w)^kw^k\\
&=&\frac{1}{n}\,[w^{n-1}]\sum_{k=0}^{\infty}\sum_{j=0}^k{n+k-1\choose k}{k\choose j}a_3^{k-j}a_6^{j}w^{j+k}\\
&=&\frac{1}{n}\sum_{k=\lfloor (n-1)/2\rfloor}^{n-1}{n+k-1\choose k}{k\choose n-k-1}a_3^{2k-n+1}a_6^{n-1-k}.
\end{eqnarray*}
\begin{Theorem}
\label{thm:family2}
Let $E$ be an elliptic curve defined over $\Q$ by the Weierstrass equation $y^2+a_3y=x^3+a_6$. Let $z=-x/y$. The $z$-power series expansion of the invariant differential $\omega_E$ of $E$ is given by
\begin{eqnarray*}
\omega_E=\sum_{n=1}^{\infty}\left(\sum_{k=\lfloor(n-1)/2\rfloor}^{n-1}{n+k-1\choose k}{k\choose n-k-1}a_3^{2k-n+1}a_6^{n-k-1}\right)z^{(3n-2)-1}\,dz.
\end{eqnarray*}
\end{Theorem}
\begin{Proof}
The invariant differential of $E$ is defined by $\displaystyle \omega_E=\frac{dw/dz}{3z^2}\,dz$, see (\ref{eq:invariant}). The result now follows as we have shown that
\[ w(z)=\sum_{n=1}^{\infty}\left(\frac{1}{n}\sum_{k=\lfloor(n-1)/2\rfloor}^{n-1}{n+k-1\choose k}{k\choose n-k-1}a_3^{2k-n+1}a_6^{n-k-1} \right)z^{3n}.\]
\end{Proof}
One can apply Atkin-Swinnerton-Dyer congruences to the elliptic curve $E$ defined over $\Q$ by the Weierstrass equation $y^2+a_3y=x^3+a_6$. In particular, if $p$ is a good prime of $E$, then using Hasse's bound, $\displaystyle |t_p|<2\sqrt{p}$, one has $t_p=0$ if $p\equiv 2$ mod $3$, moreover, $\displaystyle t_p\equiv \sum_{k=\lfloor(p-1)/6\rfloor}^{(p-1)/3}{\frac{p-1}{3}+k\choose k}{k\choose \frac{p-1}{3}-k}a_3^{2k-\frac{p-1}{3}}a_6^{\frac{p-1}{3}-k}$ mod $p$ if $p\equiv 1$ mod $3$, where $t_p=1+p-\# E(\mathbb{F}_p)$.

\begin{Remark}
Given Theorem \ref{thm:family1} and Theorem \ref{thm:family2}, one can produce a large number of congruence relations relating combinatorial objects to coefficients of modular forms, and congruences satisfied by the combinatorial quantities themselves. All one needs is applying these theorems and Atkin-Swinnerton-Dyer congruences to different families of elliptic curves described by Weierstrass equations of the form given in \S 3 and \S4.

For example, elliptic curves with $n$-torsion points are parametrized by Tate's normal form. In other words, an elliptic curve $E_n/\Q$ with a rational $n$-torsion point $(0,0)$, $n\ge 4$, is described by the following Weierstrass equation
\[y^2+(1-c)xy-by=x^3-bx^2,\;c,b\in\Z.\]
These explicit parametrizations are due to Kubert, see \cite{Kubert}. 
According to Theorem \ref{thm:family1}, the invariant differential of $E_n$ is \[\omega_{E_n}=\frac{dz}{\sqrt{1-2(1-c-bz)z+[(1-c-bz)^2+4bz]z^2}}=\sum_{n=1}^{\infty}b(n)z^{n-1}\,dz,\] where
\[b(n)=\sum_{m=0}^{n-1}\sum_{k=0}^{\lfloor m/2\rfloor}{m\choose k}{m-k\choose k}{m-2k\choose n-m-k-1}(1-c)^{2m-n-k+1}(-b)^{n-m-1}.\]
One obtains congruence relations satisfied by $b(n)$ using Corollary \ref{cor:family1}.
 \end{Remark}
 
 \hskip-18pt\emph{\bf{Acknowledgements.}}
I would like to thank the anonymous referee for his thorough reading of the article and many helpful suggestions that helped improving the manuscript.  
\bibliographystyle{plain}
\footnotesize
\bibliography{formal}
 Department of Mathematics and Actuarial Science\\ American University in Cairo\\ mmsadek@aucegypt.edu
\end{document}